\documentclass[a4paper,11pt]{amsart}
\usepackage{amsmath,amssymb,amsfonts,phonetic}
\usepackage{epsfig}
\usepackage{mathrsfs}
\newtheorem{defn}{Definition}[section]
\newtheorem{thm}{Theorem}[section]
\newtheorem{prop}{Proposition}[section]
\newtheorem{rmk}{Remark}[section]
\newtheorem{lma}{Lemma}[section]

\newtheorem{corl}{Corollary}[section]
\def\N{{\rm I\kern-0.16em N}}
\def\R{{\rm I\kern-0.16em R}}
\def\E{{\rm I\kern-0.16em E}}
\def\P{{\rm I\kern-0.16em P}}
\def\F{{\rm I\kern-0.16em F}}
\def\B{{\rm I\kern-0.16em B}}
\def\C{{\rm I\kern-0.46em C}}
\def\G{{\rm I\kern-0.50em G}}
\numberwithin{equation}{section}
\newcommand{\ud}{\mathrm{d}}
\font\eka=cmex10
\usepackage{ae}
\def\ind{\mathrel{\hbox{\rlap{%
\hbox to 7.5pt{\hrulefill}}\raise6.6pt\hbox{\eka\char'167}}}}
\begin{document}
\title[Spectral characterization of the quadratic variation]
{Spectral characterization of the quadratic variation of mixed Brownian-fractional Brownian motion}
\author[Azmoodeh and Valkeila]{Ehsan Azmoodeh \and Esko Valkeila}
\address{Department  of Mathematics and Systems Analysis\\ Aalto University\\
P.O. Box 11100, 00076 Aalto, Finland} 
\email{ehsan.azmoodeh@aalto.fi \text{ and } esko.valkeila@aalto.fi}
\thanks{ Azmoodeh is grateful to Finnish Graduate School in Stochastic and Statistic (FGSS) for financial support, 
and Valkeila acknowledges the support of Academy of Finland, grant 212875. We are grateful to an 
anonymous referee for a careful reading of the first version of this paper.}

\begin{abstract}
Dzhaparidze and Spreij \cite{ds} showed that the quadratic variation of a semimartingale 
can be approximated using a randomized periodogram. We show that the same 
approximation is valid for a special class of continuous stochastic processes. This class 
contains both semimartingales and non-semimartingales. 
The motivation comes partially from the 
recent work by Bender et al. \cite{bsv}, where it is shown that the quadratic variation of 
the log-returns determines the hedging strategy.

\medskip


\smallskip

\end{abstract}

\maketitle
\section{Introduction}
\subsection*{Spectral characterization of the bracket}
It is well-known that for a semimartingale $X$, the bracket $[X,X]$ can be identified with 
$$
[X,X]_t = \P \hbox{-} \lim _{|\pi| \to 0} \sum_{t_k\in \pi } \left( X_{t_k} - X _{t_{k-1}}\right) ^2 ,
$$
where $\pi = \left\lbrace  t_k : 0= t_0 < t_1 <\cdots < t_n = t\right\rbrace $ is 
a partition of the interval $[0,t]$, $|\pi| = \max \left\lbrace  t_k - t_{k-1} : t_k \in \pi \right\rbrace $,
and $\P \hbox{-} \lim$ means convergence in probability. Statistically speaking, the sums of squared 
increments $($\textit{realized quadratic variation}$)$ is a 
consistent estimator for the bracket as the volume of observations tends to infinity. 
Barndorff-Nielsen and Shephard \cite{bn-s} studied precision of the realized 
quadratic variation estimator for a special class of continuous semimartingales. 
They showed that sometimes the realized quadratic variation estimator can be rather 
noisy estimator. So one should seek for new estimators of the quadratic variation. \\

Dzhaparidze and Spreij \cite{ds} suggested another characterization of the bracket $[X,X]$. 
Let $\F^X$ be the filtration of $X$ and $\tau $ be a finite stopping 
time. For $\lambda \in \R $, define the \textit{periodogram} $I_\tau (X;\lambda )$ of $X$ at $\tau$ by 
\begin{equation}\label{eq:perio}
\begin{split}
I_\tau (X;\lambda ) :& = \big| \int _0^\tau e^{i\lambda s} \ud X_s \big| ^2 \\
&= 2 \ \textbf{Re} \int_{0}^{\tau} \int_{0}^{t} e^{i \lambda (t-s)} \ud X_s \ud X_t + 
[X,X]_\tau \quad (\text{ by Ito formula } ).\\
\end{split}
\end{equation}
Given $L > 0$ and $\xi$ be a symmetric random variable with a density $g_\xi$, real 
characteristic function $\varphi_\xi$, and independent of the filtration $\F^X $. 
Define the randomized periodogram by 
\begin{equation}\label{eq:p-random}
\E_\xi I_\tau \left( X; L\xi \right) = \int _\R I_\tau \left( X; Lx\right) g_\xi(x) \ud x .
\end{equation}
If the characteristic function $\varphi_\xi$ is of bounded variation, then Dzhaparidze and Spreij 
have shown  that we have the following characterization of the bracket as $L \to \infty $
\begin{equation*}\label{eq:d-s}
\E_\xi I_\tau \left( X; L\xi \right) \stackrel{\P}{\to} [X,X]_\tau  .
\end{equation*}

\subsection*{Robust Black \& Scholes pricing by hedging}
Next we give another motivation for our estimation problem. 
Bender et al. \cite{bsv} consider a class of pricing models, where the continuous stock price $S$
has the following quadratic variation as a functional of the observed path $S$:

\begin{equation}\label{eq:bsv-bracket}
\ud [S,S]_t = \sigma ^2 (S_t) \ud t.
\end{equation}

Here $\sigma : \R \to \R $ is a continuously differentiable function of  linear growth.
A typical example of this kind of stock price models is the classical Black \& Scholes model with 
constant volatility $\sigma$, where the stock price $\tilde S$ is given by 
$$
 \tilde S_t = s_0 e^{\sigma W_t -\frac12 \sigma^2 t }, 
$$
where $W$ is a standard Brownian motion. We have 
$$
\ud [\tilde S , \tilde S]_t = \sigma ^2 {\tilde S}^2 _t \ud t ,
$$
and the bracket $[\tilde S , \tilde S]$ has the form of \eqref{eq:bsv-bracket}. On the other hand, 
let $X$ be a continuous process with quadratic variation 
$[X,X] _t = \sigma ^2 t$; take for example  $X_t = \sigma W_t + \eta B^H_t $, $B^H $ is a fractional Brownian motion with Hurst parameter $H> \frac12 $,  independent of $W$, 
and $\eta $ is a constant.  Then, for the process $S_t = s_0 e^{X_t - \frac12\sigma ^2 t}$, we have that 
$$ 
\ud [S,S]_t = \sigma ^2 S^2 _t \ud t ,
$$
and again the bracket has the functional form of (\ref{eq:bsv-bracket}). Here we have examples, where the quadratic variation of the driving process $X$ determines the structure 
of the quadratic variation of the stock price. 
Moreover, if this is the case, then
Bender et al. have shown that within a fixed model class, determined by the relation \eqref{eq:bsv-bracket} the hedging of options has the same functional form 
as in the classical Black \& Scholes pricing model. The options, which can be hedged, includes
European options, path dependent options like look-back options,  and Asian options.

\subsection*{The results}
We show that the result of Dzhaparidze and Spreij holds for the mixed Brownian-fractional Brownian motion. Let  $(\Omega,\mathcal{F},\P)$ be a complete probability 
space and fix $T >0 $. Throughout the paper, we assume that $W= \{W_{t}\}_{t\in [0,T]}$ is a standard Brownian motion and  $B^{H}= \{B^{H}_{t}\}_{t \in [0,T]}$ is a 
fractional Brownian motion with Hurst parameter $H \in (\frac{1}{2},1)$, independent of the Brownian motion $W$. Define the mixed Brownian-fractional Brownian motion 
$X_{t}$ by 
\begin{equation*}
 X_{t}= W_{t}+B^{H}_{t}  \qquad t \in[0,T] .
\end{equation*}
It is known that $($see \cite{ch}$)$ the process $X$ is a $(\F ^X , \P )$ semimartingale, if $H \in( \frac34,1)$, and for $H\in \left( \frac12 , \frac34 \right]$, $X$ is 
not a semimartingale with respect to its own filtration $\F ^X$. Moreover in both cases we have
\begin{equation}\label{eq:mixed-bracket}
 [X,X]_{t}=[X]_{t}= \P \hbox{-}\lim_{ |\pi|\to 0 } \sum _{\begin{subarray}{1}
                                   t_{i}\leqq t 
                                 \end{subarray}}
 (X_{t_{i}}- X_{t_{i-1}})^{2} = t .
\end{equation}
If the partitions in (\ref{eq:mixed-bracket}) are nested, then the convergence can be strengthened to almost sure convergence. Hereafter, we always assume that the 
sequence of partitions are nested. For $\lambda\in \R $, define the complex-valued stochastic process $Y$ by 
\begin{equation*}\label{eq:y}
 Y_{t}=\int_{0}^{t}e^{i\lambda s} \ud X_{s} \qquad t\in[0,T],
\end{equation*}
where the stochastic integral is understood in a path-wise way, and it is defined by integration by parts formula $($see \cite{p}$)$:
\begin{equation*}
 \int_{0}^{t}e^{i\lambda s} \ud X_{s}=e^{i\lambda t}X_{t}-i\lambda \int_{0}^{t}X_{s}e^{i\lambda s} \ud s .
\end{equation*}
Therefore, $Y=\{Y_{t}\}_{t\in[0,T]}$ is a process with continuous sample paths. Moreover, it is straightforward to check that for $t \in [0,T]$, we have that
\begin{eqnarray*}
 [Y,Y]_t = [Y]_{t} &:= & \mbox{a.s.} \hbox{-} \lim _{|\pi|\to 0} \sum _{\begin{subarray}{1}
                                   t_{i}\leqq t  
                                 \end{subarray}}
 \Big( \left( Y_{t_k}- Y_{t_{k-1}}\right)\left(\overline{Y_{t_k} - Y_{t_{k-1}}}\right) \Big) \\ 
 &=& [\textbf{Re}\ Y]_t + [\textbf{Im}\ Y]_t = [X]_t = t,
\end{eqnarray*}
where $\overline {Y_t} $ is complex conjugate of $Y_t$ $($\cite[p.84]{p}$)$. Given $\lambda\in \R $, define the periodogram of $X$ at $T$ as $($\ref{eq:perio}$)$, i.e.
\begin{equation*}
\begin{split} 
I_{T}(X;\lambda) & = \big| \int_{0}^{T}e^{i\lambda t} \ud X_{t} \big| ^{2}\\
& = \big| e^{i\lambda T}X_{T}-i\lambda \int_{0}^{T} X_{t}e^{i\lambda t} \ud t \big|^2 \\
& =  X_T ^{2} + X_T \int_{0}^{T} i\lambda ( e^{i \lambda (T-t)} - e^{- i \lambda (T-t)} ) X_t \ud t + \lambda ^2 \big| \int_{0}^{T} e^{i \lambda t} X_t \ud t \big| ^2.
\end{split}
\end{equation*}

Assume $(\tilde{\Omega},\tilde{\mathcal{F}},\tilde{\mathbb{P}})$ be an another probability space and identify the $\sigma$-algebra 
$\mathcal{F}$ by $\mathcal{F} \otimes \{ \phi, \tilde{\Omega} \}$ on the product space 
$(\Omega \times \tilde{\Omega}, \mathcal{F} \otimes \tilde{\mathcal{F}},\mathbb{P}\otimes \tilde{\mathbb{P}})$. Let $\xi:\tilde{\Omega} \rightarrow \R$ be a real 
symmetric random variable with density $g_\xi$, and independent of the filtration $\F^X$. 
Define for any positive real number $L$ the randomized periodogram by
\begin{equation}\label{eq:r-p-m-b-fb}
 \E_{\xi} I_{T}(X;L\xi):=\int_{\mathbb{R}}I_{T}(X;Lx)g_\xi(x) \ud x.
\end{equation}

Our main result is the following.
\begin{thm}\label{t:main}
Assume that $X$ is a mixed Brownian-fractional Brownian motion, $\E_\xi I_T(X; L\xi ) $ is the randomized 
periodogram given by (\ref{eq:r-p-m-b-fb}) and 
$$
\E\xi ^2 < \infty.
$$
Then as $L \to \infty $ we have
\begin{equation*}\label{eq:main}
 \E_{\xi} I_{T}(X;L\xi) \stackrel{\P}{\to} [X,X]_T . 
\end{equation*}
\end{thm}

\begin{rmk}\label{r:d-s}
To compare our result to the results of Dzhaparidze and Spreij:\\
\begin{itemize}
 \item They take any finite stopping time $\tau$, whereas we must assume a constant stopping time $T$.
\item They assume that the characteristic function $ \varphi_\xi$ is of bounded variation, whereas instead we assume that $\xi$ is a square integrable random variable. 
\end{itemize}

Note that under our assumption of deterministic stopping time $\tau = T$, when $X$ is a Gaussian martingale, 
they can drop the condition of the bounded variation on 
$[0, \infty )$ of the characteristic function $\varphi_\xi$ of the random variable $\xi$ 
(see \cite[Remark, p.170]{ds}). 
\end{rmk}

Next we give some auxiliary material and then finish the proof.

\section{Auxiliary results}
\subsection{Path-wise Ito formula} F{\"o}llmer \cite{f1} obtained a path-wise calculus 
for continuous functions with finite quadratic 
variation. The next theorem is essentially due to F{\"o}llmer. For a nice exposition, 
and its use in finance, see Sondermann \cite{s}.

\begin{thm} \cite{s} \label{t:f-ito}
Let $X :[0,T] \rightarrow \R$ be a continuous process with continuous quadratic variation 
$[X,X]_{t}$ and $F\in C^{2}(\R)$. Then for any $t\in [0,T]$, the limit of 
the Riemann-Stieltjes sums

\begin{equation*}
 \lim_{ |\pi| \to 0} \sum _{\begin{subarray}{1}
                                   t_{i}\leqq t
                                 \end{subarray}}
F_x(X_{t_{i-1}})(X_{t_{i}}-X_{t_{i-1}}):= \int_{0}^{t}F_x(X_{s}) \ud X_{s},
\end{equation*}
exists almost surely. Moreover, we have
\begin{equation}\label{eq:f-ito}
 F(X_{t})=F(X_{0})+\int_{0}^{t}F_x(X_{s}) \ud X_{s}+\frac{1}{2} \int_{0}^{t} F_{xx}(X_{s}) \ud [X,X]_{s}.
\end{equation}
\end{thm}

\begin{lma}\label{l:mixed-f}
For the mixed Brownian-fractional Brownian motion $X$ we have 
\begin{equation}
 I_{T}(X;\lambda)=[X]_{T}+2 \ \textbf{Re} \int_{0}^{T} \int_{0}^{t}e^{i\lambda (t-s)} \ud X_{s} \ud X_{t} 
\end{equation}
where the iterated stochastic integral in the right hand side is understood in path-wise way, i.e. as the limit of 
the Riemann-Stieltjes sums. 

\end{lma}

\textbf{Proof}:
We apply Ito type formula (\ref{eq:f-ito}) to the real part $\textbf{Re}\ Y$ 
and the imaginary part $\textbf{Im}\ Y$ of the process $Y=\{Y_{t}\}_{t\in[0,T]}$ 
with the function $F(x)=x^{2}$. We obtain 
\begin{equation*}
 F( \textbf{Re}\ Y_{T} )=\int_{0}^{T} \left( 2 \int_{0}^{t} \textbf{Re}\ e^{i\lambda s} 
\ud X_{s} \right) \textbf{Re} \ e^{i\lambda t} \ud X_{t} +[\textbf{Re}\ Y]_{T}.
\end{equation*}
Similarly, we have
\begin{equation*}
  F(\textbf{Im}\ Y_{T})=\int_{0}^{T} \left( 2 \int_{0}^{t} \textbf{Re} -
\mspace{-5mu}ie^{i\lambda s} \ud X_{s} \right) 
\textbf{Re} -\mspace{-5mu}ie^{i\lambda t} \ud X_{t} +[\textbf{Im}\ Y]_{T}.
\end{equation*}
Summing the left and right hand sides of the identities, we get
\begin{equation*}
\begin{split}
\big| \int_{0}^{T}e^{i\lambda t} \ud X_{t}\big| ^{2} &
= \left( \int_{0}^{T} \textbf{Re}\ e^{i\lambda t} \ud X_{t} \right)^{2}+ \left(  \int_{0}^{T} \textbf{Re}\ 
-\mspace{-5mu}ie^{i\lambda t} \ud X_{t} \right)^{2} \\ 
 & = [X]_{T}+2\int_{0}^{T} \int_{0}^{t} \textbf{Re} \ e^{i\lambda s} \ 
\textbf{Re}\ e^{i\lambda t} \ud X_{s} \ud X_{t}\\
 & + 2\int_{0}^{T} \int_{0}^{t} -\textbf{Re}\  ie^{i\lambda s} \ -\textbf{Re}\  ie^{i\lambda t} 
\ud X_{s} \ud X_{t}\\
& = [X]_{T}+2 \ \textbf{Re} \int_{0}^{T} \int_{0}^{t} \ e^{i\lambda (t-s)} \ud X_{s} \ud X_{t}.\\
\end{split}
\end{equation*}
Note that $\textbf{Re} \ e^{i\lambda (t-s)}= \textbf{Re} \ e^{i\lambda s} \ 
\textbf{Re}\ e^{i\lambda t}+ \textbf{Re}\ -\mspace{-5mu}ie^{i\lambda s}\ 
\textbf{Re}\ -\mspace{-5mu}ie^{i\lambda t}$. 
\vskip0.25cm

\subsection{Path-wise stochastic integration in fractional Besov-type spaces}
A stochastic process $X$ is a semimartingale if and only if one has a version of the Lebesgue 
dominated convergence theorem $($see \cite{p}$)$. Fractional 
Brownian motion is not a semimartingale, and hence the stochastic integral with 
respect to fractional Brownian motion $B^H$ must be defined. Using the smoothness of 
 the sample paths of the fractional Brownian motion $B^H$, when $H \in (\frac{1}{2},1)$, 
one can define the so-called \textit{generalized Lebesgue-Stieltjes integral}. 
For more information, see \cite{nr}, \cite{z} and \cite{m}.

\begin{defn}\cite{nr}
Fix $0<\alpha< \frac{1}{2}$.\\
(i) For $f:[0,T]\to \mathbb{R}$, define
\begin{equation*}
 \Vert f \Vert_{\alpha,1}:= \int_{0}^{T} \frac{|f(t)|}{t^\alpha} \ud t + 
\int_{0}^{T} \int_{0}^{t} \frac{|f(t)-f(s)|}{(t-s)^{\alpha + 1}} \ud s \ud t,
\end{equation*}
and
\begin{equation*}
 W^{\alpha,1}_{0}[0,T]=\{f:[0,T]\to \mathbb{R} \ ; \ \Vert f \Vert_{\alpha,1}< \infty \}.
\end{equation*}
(ii) Also, for $f:[0,T]\to \mathbb{R}$, define

\begin{equation*}
 \Vert f \Vert_{1-\alpha, \infty,T}:= 
\sup _{0<s<t<T}\left(\frac{\vert f(t)-f(s)\vert}{(t-s)^{1-\alpha}}+\int_{s}^{t} 
\frac{\vert f(y)-f(s)\vert}{(y-s)^{2-\alpha}}\ud y \right),
\end{equation*}
and
\begin{equation*}
 W^{1-\alpha,\infty}_{T}[0,T]:=\{f:[0,T]\to \mathbb{R} \ ; \ \Vert f \Vert_{1-\alpha, \infty,T}< \infty \}.
\end{equation*}
\end{defn}

Denote by $C^{\lambda}[0,T]$ the space of $\lambda$-H\"older continuous functions on the interval $[0,T]$. Then $\forall \epsilon > 0$, we have the inclusions
\begin{equation*}
\begin{split} 
& C^{1-\alpha+\epsilon}[0,T]  \subseteq  W^{1-\alpha,\infty}_{T}[0,T]  \subseteq C^{1-\alpha}[0,T] \\
                            & C^{\alpha + \epsilon} [0,T] \subseteq W^{\alpha,1}_{0} [0,T].\\
\end{split}
\end{equation*}
Recall that almost surely sample paths of $B^H$ for any $0<\gamma < H$, belong to $ C^{\gamma }[0,T] $. This follows from the Kolmogorov continuity theorem. Hence the
sample paths of $B^H$ belong to $ W^{\alpha,\infty}_{T}[0,T]$ for any $ 0 < \alpha < H $ .\\

In the following  $D^{\alpha}_{t^{-}} (\mbox{resp.} D^{\alpha}_{0^+})$ stand for right-sided (resp. left-sided) 
fractional derivatives (\cite{skm}). For $g\in W^{1-\alpha,\infty}_{T}[0,T]$, define
\begin{equation*}
 \Lambda_{\alpha}(g):=\frac{1}{\Gamma(1-\alpha)}\sup _{0<s<t<T} 
\vert (D^{1-\alpha}_{t^{-}}g_{t^{-}})(s) \vert \\
\le \frac{1}{\Gamma(1-\alpha)\Gamma(\alpha)} \Vert g \Vert_{1-\alpha, \infty,T}.
\end{equation*}

\begin{defn}\cite{nr}
Fix  $ 0 <  \alpha < \frac{1}{2}$. Let $f\in  W^{\alpha,1}_{0}[0,T] $ and $ 
g\in W^{1-\alpha,\infty}_{T}[0,T]$. Then the Lebesgue integral 
\begin{equation*}
\int_{0}^{T} D^{\alpha}_{0^{+}} f_{0^{+}}(t) D^{1- \alpha}_{T^{-}}g_{T^{-}}(t) \ud t
\end{equation*}
exists, and we can define the generalized Lebesgue-Stieltjes integral by
\begin{equation*}
\int_{0}^{T} f_{t}dg_{t} := \int_{0}^{T} D^{\alpha}_{0^{+}} f_{0^{+}}(t) D^{1- \alpha}_{T^{-}}g_{T^{-}}(t) \ud t,
\end{equation*}
where $f_{0^{+}}(t)= f(t) - f(0^+) $ and $ g_{T^{-}}(t)=g(T^-) - g(t)$.
\end{defn}

\begin{rmk} \cite{m}, \cite{z}
The definition of the generalized Lebesgue-Stieltjes integral does not depend on the choice of $\alpha$. 
\end{rmk}

\begin{rmk}\cite{z}\label{r:coincide}
If $f$ and $ g$ are H\"older continuous of orders $\alpha$ and $\beta$ with $\alpha + \beta > 1$, 
then the generalized Lebesgue-Stieltjes integral exists and coincides 
with the Riemann-Stieltjes integral. This fact is based on the integration theory developed by Young \cite{y}.
\end{rmk}
Since fractional Brownian motion is not a semimartingale, the next theorem and corollary can be used 
instead of the Lebesgue dominated convergence theorem for fractional Brownian motion.
\begin{thm}\cite{nr}
 Let $g\in W^{1-\alpha,\infty}_{T}[0,T]$ and $f\in  W^{\alpha,1}_{0}[0,T]$. Then we have the estimate
\begin{equation*}
 \Big| \int_{0}^{T}f_{t} \ud g_{t} \Big| \le \Lambda_{\alpha}(g) C_{\alpha,T}\Vert f \Vert_{\alpha,1}
\end{equation*}
for some constant $C=C_{\alpha,T}$.
\end{thm}

\begin{corl}\label{c:n-r}
Assume $ f, f^n  \in  W^{\alpha,1}_{0}[0,T]$, and $\Vert f^n-f \Vert _{\alpha,1} \to 0$ 
as $n\to \infty$ for some $\alpha \in (1-H,\frac{1}{2})$. Then as $n \to \infty$
\begin{equation*}
\int_{0}^{T} f^{n}_{t} \ud B^{H}_{t}  \to   \int_{0}^{T} f_t  \ud B^{H}_{t}, \quad \text{a.s.}
\end{equation*}
\end{corl}

Next we use this machinery to prove a stochastic Fubini type result that is 
cornerstone of the proof of our main result.

\begin{prop}\label{l:mixed-fubini}
Assume that $X = W + B^H $ is a mixed Brownian-fractional Brownian motion, 
and $\xi$ is a square integrable random variable with a density $g_\xi$, and independent 
of the filtration $ \F^X$. Then the iterated integrals
\begin{gather*}
I_{1}=\int_{\mathbb{R}}\left(\int_{0}^{T}\phi_{W}(t,x)dB^{H}_{t}\right)g_\xi(x)dx \notag \\
I_{2}=\int_{0}^{T} \left(\int_{\mathbb{R}}\phi_{W}(t,x)g_\xi(x)dx\right)dB^{H}_{t}
\end{gather*}
exist, and moreover $ I_{1}= I_{2}$ almost surely, where $ \phi_{W}(t,x)= \int _{0}^{t} e^{ix(t-s)} dW_s $ 
and the stochastic integrals in $I_1$ and $I_2$ are 
understood in path-wise way as the limit of Riemann-Stieltjes sums.
\end{prop}

\textbf{Proof}: We split the proof into four steps.\\

\textit{Step $1$: The existence of $ I_1$}. Using integration by parts formula and simple 
manipulations, we see that the sample paths of the complex-valued stochastic 
process $\{ \phi_{W}(t,x)\} _{t\in[0,T]}$ parametrized by $x \in \mathbb{R}$, 
are  H\"older continuous of any order less than half almost surely. Moreover, for any 
$\alpha \in (0,\frac{1}{2})$, $x\in \mathbb{R},$ and $ s,t \in [0,T]$, we have 
\begin{equation}
\begin{split}
\big| \phi_{W}(t,x) - \phi_{W}(s,x) \big| &\le  C(\omega,T)  (1+|x|+|x|^2) 
\big| t - s \big|^{\alpha} \text{ and }\\
 \big| \phi_{W}(t,x) \big| & \le C(\omega,T) ( 1 + T |x| ), \label{eq:existence}\\
\end{split}
\end{equation}
where $ C(\omega,T) $ is an almost surely finite and positive random variable that may 
be different from line to line. Hence the interior stochastic integral in 
$I_{1}$ can be defined as limit of Riemann-Stieltjes sums almost surely by using the Young integration theory 
(see Remark \ref{r:coincide}). Note that 
by \eqref{eq:existence}, the function 
\begin{equation*}
 x \mapsto \int_{0}^{T} \phi_{W}(t,x) \ud B^{H}_t 
\end{equation*}
is integrable with respect to the measure $g_{\xi}(x) \ud x$ almost surely.\\

\textit{Step $2$: The existence of $ I_2$}.
Using $($\ref{eq:existence}$)$, we see that for any $\alpha \in (0,\frac{1}{2}) $,
\begin{equation*}
\Big| \int_{\mathbb{R}} \phi_{W}(t,x) g_{\xi}(x) \ud x - 
\int_{\mathbb{R}} \phi_{W}(s,x) g_{\xi}(x) \ud x \Big| 
\le  C(\omega,T) ( 1 + 2 \E \xi^2) \big| t-s \big|^\alpha.
\end{equation*}
Therefore, the stochastic integral $I_2$ can be defined as limit of 
the Riemann-Stieltjes sums almost surely.\\

\textit{Step $3$}. Define for any $N \in \mathbb{N}$,
\begin{equation*}
\begin{split} 
I^{N}_{1} & =\int_{-N}^{N}\left(\int_{0}^{T}\phi_{W}(t,x) \ud B^{H}_{t}\right)g_\xi (x) \ud x  \\
I^{N}_{2} & =\int_{0}^{T}\left(\int_{-N}^{N}\phi_{W}(t,x)g_\xi(x) \ud x\right) \ud B^{H}_{t}.
\end{split}
\end{equation*}
Clearly $I^{N}_{1}$ converges to $I_1$ almost surely as $N$ tends to infinity. 
We aim to show that  $I^{N}_{1}=I^{N}_{2}$ almost surely. By definition of the Riemann 
integral, there exists a sequence of partitions 
$\{ \pi ^{N} _n \}_{n=1}^{\infty}$ of the interval $[-N,N]$, 
such that $| \pi ^{N}_n | \to 0 \text{ as } n \to \infty$ and 
\begin{equation*}
 I^{N}_1 = \lim _{n \to \infty} \sum _{x^{n}_i \in \pi ^{N}_n} 
\left( \int _{0}^{T} \phi _{W} (t,x^{n}_{i-1}) \ud B^{H}_t \right) g_\xi (x^{n}_{i-1}) \Delta x^{n}_{i}, 
\end{equation*}
and
\begin{equation*}
\lim _{n \to \infty} \sum _{x^{n}_i \in \pi ^{N}_n}
\phi _{W} (t,x^{n}_{i-1}) g_\xi (x^{n}_{i-1}) \Delta x^{n}_{i} = \int_{N}^{N} \phi_{W} (t,x)g_\xi (x) \ud x
\end{equation*}
hold. Assume that $\pi ^{N}_{n} =\{ -N = x^{n}_{0}  < x^{n}_{1}< \cdots < x^{n}_{k_n} = N \}$. 
For each  $x^{n}_{i-1} \in \pi ^{N}_{n},\  1 \le i \le k_{n} + 1$, we 
can find a sequence $\{ \pi ^{T,x^{n}_{i-1}} _{m} \}$ of partitions of 
the interval $[0,T]$, such that $| \pi ^{T,x^{n}_{i-1}} _{m} | \to 0 \text{ as } m \to \infty $ and
\begin{equation*}
\lim _{m \to \infty} \sum _{t^{m}_{j} \in \pi ^{T,x^{n}_{i-1}} _{m}} 
\phi _{W}(t^{m}_{j-1}, x^{n}_{i-1}) \Delta B^{H}_{t^{m}_{j}} = 
\int _{0}^{T} \phi _{W}(t, x^{n}_{i-1}) \ud B^{H}_{t}.
\end{equation*}
On the other hand, there is another sequence $ \{ \hat{\pi}_m ^{T} \}$ of partitions of 
the interval $[0,T]$ such that
\begin{equation*}
\begin{split} 
\lim _{m \to \infty} \sum _{\hat{t}^{m}_{j} \in \hat{\pi}_{m}^{T}�} 
\left( \int_{-N}^{N} \phi_{W}(\hat{t}^{m}_{j-1},x) g_\xi (x) \ud x \right) & \Delta B^{H}_{\hat{t}^{m}_{j}}\\ 
&= \int_{0}^{T}  \left( \int_{-N}^{N}  \phi_{W}(t,x)g_\xi(x) \ud x \right) \ud B^{H}_{t}.\\
\end{split}
\end{equation*}
Let 
\begin{equation*}
 \pi ^{T,n}_{m} = \cup _{i=1}^{k_{n}+1} \pi ^{T,x^{n}_{i-1}} _{m} \quad 
\text{ and }\quad  \pi^{T}_{m} = \hat{\pi}_m ^{T} \cup \pi ^{T,n}_{m}.
\end{equation*}

Therefore, for any $n \in \mathbb{N}$, the partition $\pi^{T}_{m} $ of the interval $[0,T]$ 
contains all points of the partitions $\pi ^{T,n}_{m} $ and 
$\hat{\pi}_m ^{T}$, and denote the points of $ \pi^{T}_{m}$ by $t^{m}_{k}, \ k= 0, \cdots, l_m$. Then  for any $x^{n}_{i-1} \in \pi ^{N}_{n} $, we can write
\begin{equation*}
 \lim _{m \to \infty} \sum _{t^{m}_{j} \in \pi^{T}_{m}} \phi _{W}(t^{m}_{j-1}, x^{n}_{i-1}) 
\Delta B^{H}_{t^{m}_{j}} = \int _{0}^{T} \phi _{W}(t, x^{n}_{i-1}) \ud B^{H}_{t}.
\end{equation*}
Now for $n,m \in \mathbb{N}$, we can have the estimate 
\begin{equation*}
| I^{N}_1- I^{N}_2 | \le | I^{N}_1 - \Delta _{n,m}| + | I^{N}_2 - \Delta _{n,m}| := 
A_{n,m} + B_{n,m},
\end{equation*}
where 
\begin{equation*}
 \Delta _{n,m}:=\sum _{x^{n}_i \in \pi ^{N}_{n}}  \sum _{t^{m}_{j} \in \pi^{T}_{m}} 
\phi _{W}(t^{m}_{j-1}, x^{n}_{i-1}) \Delta B^{H}_{t^{m}_{j}} g_\xi (x^{n}_{i-1}) \Delta x^{n}_{i}.
\end{equation*}

For the first term $A_{n,m}$, we have 
\begin{equation*}
\begin{split}
| A_{n,m} | & \le \Big| I^{N}_1 - \sum _{x^{n}_i \in \pi ^{N}_n} 
\big( \int _{0}^{T}\phi _{W} (t,x^{n}_{i-1}) \ud B^{H}_t \big) g_\xi (x^{n}_{i-1}) \Delta x^{n}_{i} \Big| \\
& + \sum _{x^{n}_i \in \pi ^{N}_n} \Big| \int _{0}^{T}\phi _{W} (t,x^{n}_{i-1}) \ud B^{H}_t\\
 & - \sum _{t^{m}_{j} \in \pi^{T}_{m}} \phi _{W}(t^{m}_{j-1}, x^{n}_{i-1}) 
\Delta B^{H}_{t^{m}_{j}} \Big| g_\xi (x^{n}_{i-1}) \Delta x^{n}_{i}.\\
\end{split}
\end{equation*}

For fix $n \text{ and } x^{n}_{i-1}$, when $ m$ tends to infinity, we have
\begin{equation*}
\Big| \int _{0}^{T}\phi _{W} (t,x^{n}_{i-1}) \ud B^{H}_t - 
\sum _{t^{m}_{j} \in \pi^{T}_{m}} \phi _{W}(t^{m}_{j-1}, x^{n}_{i-1}) \Delta B^{H}_{t^{m}_{j}} \Big| \to 0.
\end{equation*}
Therefore, $\lim_{n \to \infty} \lim_{m \to \infty}  A_{n,m} = 0$. Similarly, for the second term $B_{n,m}$
\begin{equation*}
\begin{split}
| B_{n,m} | & \le \Big| I^{N}_2 - \sum _{t^{m}_{j} \in \pi^{T}_{m}} 
\left( \int_{-N}^{N} \phi_W (t^{m}_{j-1},x)g_\xi (x) \ud x \right ) \Delta B^{H}_{t^{m}_{j}} \Big| \\
& + \Big| \Delta _{n,m}- \sum _{t^{m}_{j} \in \pi^{T}_{m}} 
\left( \int_{-N}^{N} \phi_{W} (t^{m}_{j-1},x)g_\xi (x) \ud x \right) \Delta B^{H}_{t^{m}_{j}} \Big|.\\
\end{split}
\end{equation*}
So, it is enough to show that the second term in the right hand side converges to $0$ as $n,m$ tend to 
infinity. Note that the second term can be written 
as 
\begin{multline*}
\Big| \Delta _{n,m}-\sum _{t^{m}_{j} \in \pi^{T}_{m}} 
\left( \int_{-N}^{N} \phi_{W} (t^{m}_{j-1},x)g_\xi (x) \ud x \right) \Delta B^{H}_{t^{m}_{j}}\Big| 
=\Big|\sum _{i=0}^{k_n} \int_{x^{n}_{i-1}}^{x^{n}_{i}} f_m (x,x^{n}_{i-1}) \ud x \Big|,
\end{multline*}
where 
\begin{equation*}
 f_m (x,x^{n}_{i-1}) = \sum_{t^{m}_j \in \pi^{T}_{m}} 
\big( \phi _{W}(t^{m}_{j-1}, x^{n}_{i-1})g_\xi (x^{n}_{i-1}) -
\phi _{W}(t^{m}_{j-1}, x)g_\xi(x) \big) \Delta B^{H}_{t^{m}_{j}}. 
\end{equation*}
So, when $m$ tends to infinity, we have that
\begin{equation*}
 f_m (x,x^{n}_{i-1})  \to \int_{0}^{T} 
\Big( \phi_W (t, x^{n}_{i-1})g_\xi (x^{n}_{i-1})-\phi _W (t,x) g_\xi(x) \Big) \ud B^{H}_t.
\end{equation*}

Moreover, for each $ 1 \le i \le k_n$, the sequence $f_m (x,x^{n}_{i-1}) $ has an 
integrable dominant with respect to variable $x$. To see this, take 
$\theta \in ( \frac{1}{2}, H) $ and $ \lambda \in (0,\frac{1}{2} ) $ such 
that $ \theta + \lambda = 1 + \epsilon$. Then
\begin{equation*}
\begin{split}
\Big| f_m (x,x^{n}_{i-1}) \Big| & \le \Big| \sum_{t^{m}_j \in \pi^{T}_{m}} 
\big( \phi _{W}(t^{m}_{j-1}, x^{n}_{i-1})g_\xi (x^{n}_{i-1}) -\phi _{W}(t^{m}_{j-1}, x)g_\xi(x) \big)  
\Delta B^{H}_{t^{m}_{j}} \\
& - \int_{0}^{T} \big( \phi_W (t, x^{n}_{i-1})g_\xi (x^{n}_{i-1}) - \phi _W (t,x) g_\xi(x) \big) 
\ud B^{H}_t \Big| \\
& + \Big| \int_{0}^{T} \big( \phi_W (t, x^{n}_{i-1})g_\xi (x^{n}_{i-1}) - 
\phi _W (t,x) g_\xi(x) \big) \ud B^{H}_t \Big| \\
& \le C |\pi^{T}_{m} |^{\epsilon} \Vert B^{H} \Vert_{C^{\theta}[0,T]} 
\Vert \phi_W (t, x^{n}_{i-1})g_\xi (x^{n}_{i-1})-\phi _W (t,x) g_\xi(x) \Vert_{C^{\lambda}[0,T]}\\
& \le C |\pi^{T}_{m} |^{\epsilon} \Vert B^{H} \Vert_{C^{\theta}[0,T]} 
\Big[ \Vert \phi_W (t, x^{n}_{i-1})g_\xi (x^{n}_{i-1}) \Vert_{ C^{\lambda}[0,T]}\\
& + \Vert \phi_W (t, x)g_\xi (x) \Vert_{ C^{\lambda}[0,T]} \Big].\\
\end{split}
\end{equation*}
By the inequalities were obtained in $($\ref{eq:existence}$)$, we see that 
\begin{equation*}
\Vert \phi_W (t, x)g_\xi (x) \Vert_{ C^{\lambda}[0,T]} \le C(\omega,T) 
( 1 + |x| + |x|^2 ) g_\xi (x) \in L^{1} [-N,N].
\end{equation*}
 Therefore, by the Lebesgue dominated convergence theorem, we have that as $m$ tends to infinity
\begin{equation*}
\begin{split}
\int_{x^{n}_{i-1}}^{x^{n}_{i}} f_m (x,x^{n}_{i-1})&  \ud x \\
 & \to  \int_{x^{n}_{i-1}}^{x^{n}_{i}} \left( \int_{0}^{T} \big( \phi_W (t, x^{n}_{i-1})g_\xi (x^{n}_{i-1}) 
- \phi _W (t,x) g_\xi(x) \big) \ud B^{H}_t \right) \ud x.\\
\end{split}
\end{equation*}
Therefore, as $n$ tends to infinity, we have 
\begin{equation*}
\sum _{i=0}^{k_n} \int_{x^{n}_{i-1}}^{x^{n}_{i}} \left( \int_{0}^{T} 
\big( \phi_W (t, x^{n}_{i-1})g_\xi (x^{n}_{i-1}) - \phi _W (t,x) g_\xi(x) \big) \ud B^{H}_t \right) \ud x \to 0.
\end{equation*}
Hence, we have shown that $\lim_{n \to \infty} \lim_{m \to \infty} B_{n,m} = 0$.
\begin{rmk}
The result of this step can be derived from theorem 2.6.5 of \cite[p.177]{m}, with some modifications.
\end{rmk}
\vskip0.10cm

\textit{Step $4$}: We want to show that $I^{N}_{2}$ converges to $ I_{2}$ as $N$ tends to 
infinity. Clearly, the difference is
\begin{equation}\label{eq:step4}
 \big| I_{2}-I^{N}_{2} \big|  = \big| \int_{0}^{T} u_{t}^{N} \ud B^{H}_{t} \big|
\end{equation}
where 
\begin{equation*}
u^{N}_{t}:=\int_{[-N,N]^{c}}\phi_{W}(t,x)g_\xi(x) \ud x.
\end{equation*}
According to Corollary \ref{c:n-r}, it is sufficient to show that for some 
$\alpha \in (1- H , \frac{1}{2})$, the sequence $u^{N} \in W^{\alpha,1}_{0} [0,T]$ and 
$\Vert u^{N} \Vert _{\alpha,1} \to 0 $. Note that the sample paths of the process 
$u^N$ are H\"older continuous of any order less than half almost surely. Therefore, by 
the Remark \ref{r:coincide}, the stochastic integral appears in \eqref{eq:step4} coincides 
with the Riemann-Stieltjes integral. Now, for any $ \alpha \in (1- H , \frac{1}{2})$, using \eqref{eq:existence} 
and the assumption
$\E \xi^2 < \infty$, we have
\begin{equation*}
\int_{0}^{T} \frac{|u^N_{t}|}{t^\alpha} \ud t \to 0 \quad \text{as} \quad  N \to \infty,
\end{equation*}
by the Lebesgue dominated convergence theorem. For the second term, we take a positive 
real number $ \beta \in ( \alpha,\frac{1}{2})$. Then using \eqref{eq:existence}, 
we have 
\begin{equation*}
\begin{split}
\int_{0}^{T} \int_{0}^{t} \frac{|u^N_{t} - u^N_{s} |}{(t-s)^{\alpha + 1}} \ud s 
\ud t  &\le C(\omega,T) \int_{0}^{T} \int_{0}^{s} \frac{1}{(t-s)^{\alpha + 1 - \beta}} \ud s \ud t \\
& \int_{\mathbb{R}} \textbf{1}_{[-N,N]^c} ( 1 +|x| + |x|^2 ) g_\xi (x) \ud x\\
& \to 0 \quad \text{as} \quad N \to \infty,
\end{split}
\end{equation*}
by the Lebesgue dominated convergence theorem, since $\alpha + 1 - \beta < 1 $. 
Hence, we have shown that $\Vert u^N \Vert_{\alpha,1} \to 0 $ as $N$ tends to infinity.\\

\section{Proof of the Main Result}

Let $\varphi _\xi$ stands for the real valued characteristic function of $\xi$. 
Then the parametrized stochastic Fubini theorem for semimartingales $($see \cite{p}$)$, Lemma $($\ref{l:mixed-f}$)$ and Proposition $($\ref{l:mixed-fubini}$)$ allow us to write the randomized 
periodogram of the mixed Brownian-fractional Brownian motion $X$ as



\begin{equation*}
\begin{split}
\E_{\xi} & I_{T}(X;L\xi) = 2 \int_{0}^{T} \int_{0}^{t}\varphi _\xi (L(t-s)) \ud 
X_{s} \ud X_{t} +  [X]_{T} \\
 &= 2 \int_{0}^{T} \int_{0}^{t}\varphi_\xi(L(t-s)) \ud W_{s} \ud W_{t}+ 
2 \int_{0}^{T} \int_{0}^{t}\varphi _\xi(L(t-s)) \ud W_{s} \ud B^{H}_{t}\\
 & + 2 \int_{0}^{T} \int_{0}^{t}\varphi _\xi(L(t-s)) \ud B^{H}_{s} \ud 
W_{t}+ 2 \int_{0}^{T} \int_{0}^{t}\varphi _\xi(L(t-s)) \ud B^{H}_{s} \ud B^{H}_{t} + [X]_{T}\\
 & = 2J_{1}+2 J_{2}+ 2 J_{3}+ 2J_{4} + [X]_{T}. \\
\end{split}
\end{equation*}
Next, we show that as $L \to \infty$

\begin{equation*}
 J_k \stackrel{\P}{\to} 0, \quad k=1,2,3,4,
\end{equation*}
using the facts that 
\begin{equation*}
|\varphi_{\xi}| \le 1 \quad \text{and} \quad \varphi_{\xi}(L(t-s)) \to 0 \quad 
\text{for}\quad s < t \quad \text{as} \quad L \to \infty.
\end{equation*}
\vskip0.10cm

$J_1 \stackrel{\P}{\to} 0 $ :\\
By Ito isometry, we have 
\begin{equation*}
 \E J^{2}_{1} = \int _{0}^{T}\int_{0}^{t}\varphi ^{2}_\xi(L(t-s))\ud s \ud t \to 0 
\quad \text{as} \quad L \to \infty.
\end{equation*}

\vskip0.10cm

$J_2 \stackrel{\P}{\to} 0 $ :\\
Since Brownian motion $W$ and fractional Brownian motion $B^H$ are independent, we can compute 
\begin{equation*}
\begin{split}
 \E  \Bigg( & \int_{0}^{T} \int_{0}^{t}\varphi _\xi(L(t-s)) \ud W_{s} \ud B^{H}_{t} \Bigg)^{2}\\
&= \E \Bigg( \E \Big( \int_{0}^{T} \int_{0}^{t}\varphi _\xi(L(t-s)) \ud W_{s} \ud B^{H}_{t}\Big) ^{2} \ \Big| \ \F^{W}_T \Bigg)\\
&= H(2H-1) \E \int_{0}^{T} \int_{0}^{T} \vert u-v \vert ^{2H-2} \\
& \hspace{3.5cm} \int_{0}^{u}\varphi _\xi(L(u-s)) \ud W_{s} \int_{0}^{v}\varphi _\xi(L(v-s))\ud W_{s} \ud u \ud v \\
&= H(2H-1)\int_{0}^{T} \int_{0}^{T} \vert u-v \vert ^{2H-2} \int_{0}^{u \wedge v}\varphi _\xi(L(u-s)) \varphi _\xi(L(v-s))\ud s  \ud u \ud v\\
& \to 0 \quad \text{as} \quad L \to \infty,\\
\end{split}
\end{equation*}
by the Lebesgue dominated convergence theorem.\\

$J_3 \stackrel{\P}{\to} 0 $ :\\
Similar to the case $J_2$, we can compute 
\begin{equation*}
\begin{split}
\E \Bigg( & \int_{0}^{T}\int_{0}^{t} \varphi_{\xi}(L(t-s)) \ud B^{H}_{s} \ud W_{t}\Bigg)^{2}\\
&=\E \Bigg( \E \Big(  \int_{0}^{T}\int_{0}^{t} \varphi_{\xi}(L(t-s)) \ud B^{H}_{s}\ud W_{t} \Big)^{2} \ \Big| \ \F^{B^H}_{T} \Bigg)\\
& = H(2H-1) \int_{0}^{T}\int_{0}^{t}\int_{0}^{t}\varphi_{\xi}(L(t-u))\varphi_{\xi}(L(t-v))|u-v|^{2H-2}\ud u \ud v \ud t\\
 & \to 0 \quad \text{as} \quad L \to \infty.\\
\end{split}
\end{equation*}
\vskip0.10cm

$J_4 \stackrel{\P}{\to} 0 $ :\\
By theorem $4.1$, \cite{dk} and the Lebesgue dominated convergence theorem, we have 
\begin{equation*}
\begin{split}
\E \Bigg( & \int_{0}^{T}\int_{0}^{t} \varphi_{\xi}(L(t-s)) \ud B^{H}_{s} \ud B^{H}_{t}\Bigg)^{2}\\
 & = (H(2H-1))^{2} \int_{0}^{T}\int_{0}^{T}\int_{0}^{v}\int_{0}^{u} \varphi_{\xi}(L(u-s))\varphi_{\xi}(L(v-t))\\
& \hspace{3.5cm} |t-s|^{2H-2}|u-v|^{2H-2}\ud s \ud t \ud u \ud v\\
& \to 0 \quad \text{as} \quad L \to \infty.\\
\end{split}
\end{equation*}
\section{More properties and remarks}

Assume that $X$ is a mixed Brownian-fractional Brownian motion, i.e. $X_{t}= W_{t}+B^{H}_{t}$. Let $ \pi =\{t_{0}=0 < t_{1}< \cdots <t_{n}=T\}$ be a partition of the 
interval $[0,T]$. Then, we have the following properties of the realized quadratic variation estimator.
\begin{itemize}
 \item  Using the Ito type formula $($\ref{eq:f-ito}$)$, we have a representation for the error term, denoted by $e^1$, as
\begin{equation*}
 \sum _{t_{k}\in \pi}(X_{t_{k}} - X_{t_{k-1}})^{2}- [X,X]_T = 2 \sum _{t_{k}\in \pi} \int _{t_{k-1}}^{t_{k}}\int _{t_{k-1}}^{t} \ud X_{s} \ud X_{t}.
\end{equation*}
\item Hence, for the error term $e^1$ of the realized quadratic variation estimator, we obtain
\begin{equation*}
\begin{split}
\E (e^1) &= \E \Big(  2 \sum _{t_{k}\in \pi} \int _{t_{k-1}}^{t_{k}}\int _{t_{k-1}}^{t} \ud X_{s} \ud X_{t} \Big) \\
&= \E \Big( 2 \sum _{t_{k}\in \pi} \int _{t_{k-1}}^{t_{k}}\int _{t_{k-1}}^{t} \ud B^{H}_{s} \ud B^{H}_{t} \Big) = \sum _{t_{k}\in \pi} (\Delta t_{k})^{2H}.
\end{split}
\end{equation*}
This implies that the realized quadratic variation is a \textit{biased} estimator of the quadratic variation $[X,X]$.
\item Moreover, its variance is given by
\begin{equation*}
\begin{split} 
\mathbb{V}\text{ar} (e^1) & = \mathbb{V}\text{ar} \Big( 2 \sum _{t_{k}\in \pi} \int _{t_{k-1}}^{t_{k}}\int _{t_{k-1}}^{t} \ud X_{s} \ud X_{t}\Big)\\
&= \sum _{k=1}^{n} \mathbb{V}\text{ar} \Big( 2 \int _{t_{k-1}}^{t_{k}}\int _{t_{k-1}}^{t} \ud X_{s} \ud X_{t}\Big)\\
& + 2 \sum _{\begin{subarray}{1}
                                   1 \le i,j \le n\\ i<j  
                                 \end{subarray}}
\mathbb{C}\text{ov} \Big(  2 \int _{t_{i-1}}^{t_{i}}\int _{t_{i-1}}^{t} \ud X_{s} \ud X_{t}, 2 \int _{t_{j-1}}^{t_{j}}\int _{t_{j-1}}^{t} \ud X_{s} \ud X_{t}\Big)\\
&=\sum _{k=1}^{n} 2 \Big( (\Delta t_{k})+ (\Delta t_{k})^{2H} \Big)^{2}+ 4 \sum _{\begin{subarray}{1}
                                                                                    1 \le i,j \le n\\ i<j
                                                                                  \end{subarray}}
\Big( \E(\Delta B^{H}_{t_{i}})^2 (\Delta B^{H}_{t_{j}})^2 \Big)^{2}\\
&=\sum _{k=1}^{n} 2 \Big( (\Delta t_{k})+ (\Delta t_{k})^{2H} \Big)^{2}\\
&+\sum _{\begin{subarray}{1}
         1 \le i,j \le n\\ i<j   
          \end{subarray}}
\Big( (t_{j}- t_{i-1})^{2H}+(t_{j-1}- t_{i})^{2H} - (t_{j}- t_{i})^{2H}- (t_{j-1}- t_{i-1})^{2H}\Big)^{2}.
\end{split}
\end{equation*}
\item For the special case of equidistant partition $\pi_n = \{ \frac{kT}{n}; k=0,1,\cdots,n \} $, the mean and the variance of the error term $e^1=e_{\pi_n}^1$ 
take the forms 
\begin{equation*}
\begin{split} 
\E (e_{\pi_{n}}^1) &= T^{2H}n^{1-2H},\\
\mathbb{V}\text{ar} (e_{\pi_{n}}^1) &= 2n \Big( (\frac{T}{n}) + (\frac{T}{n})^{2H} \Big)^{2}\\
&+ (\frac{T}{n})^{4H} \sum _{\begin{subarray}{1}
                           1 \le i,j \le n\\ i<j  
                                 \end{subarray}}
\Big( (j-i-1)^{2H}+(j-i+1)^{2H}-2(j-i)^{2H} \Big)^{2}\\
\end{split}
\end{equation*}
Therefore, we have the asymptotic behaviors 
\begin{eqnarray*}
\E (e_{\pi_{n}}^1)  & \sim & T \quad \text{as} \quad H \downarrow {\frac{1}{2}},\\
\E (e_{\pi_{n}}^1)  &\to & 0 \quad \text{as} \quad  n \to \infty \quad \forall \ H >\frac{1}{2},\\
\mathbb{V}\text{ar} (e_{\pi_{n}}^1) &\sim & 2n \Big( (\frac{T}{n}) \Big)^{2} = 8 \frac{T^{2}}{n} \quad \text{as} \quad  H \downarrow {\frac{1}{2}}.
\end{eqnarray*}
Hence, we see that
\begin{equation*}
  \lim _{n \to \infty} \lim _{H \downarrow {\frac{1}{2}}}  \mathbb{V}\text{ar} (e_{\pi_{n}}^1) = 0, 
\end{equation*}
whereas for two independent Brownian motions $ W^{1}$ and $ W^{2}$, with  $Z_{t}= W^{1}_{t}+ W^{2}_{t}$ and a simple computation we have 
\begin{equation*}
  \mathbb{V}\text{ar}  \Big( 2 \int _{0}^{T} \int_{0}^{t} \ud Z_{s} \ud Z_{t}\Big) = \mathbb{V}\text{ar} \Big( Z^{2}_{T} - [Z,Z]_{T} \Big) = 8T^{2}.
\end{equation*} 
\end{itemize}

For randomized periodogram, we have the following properties.
\begin{itemize} 
\item The error term, denoted by $e^2$, of the randomized periodogram takes a form as
\begin{equation*}
\E_\xi I_T (X;L\xi) - [X,X]_T = 2 \int_{0}^{T} \int_{0}^{t} \varphi_\xi (L(t-s)) \ud X_s \ud X_t.
\end{equation*}
\item The mean of the error term $e^2$ can be computed as
\begin{equation*}
\begin{split}
\E (e^2) & = \E \big( 2 \int_{0}^{T} \int_{0}^{t} \varphi_\xi (L(t-s)) \ud X_s \ud X_t \big) \\
 & = 2H (2H-1) \int_{0}^{T} \int_{0}^{t} \varphi_\xi (L(t-s)) |t-s|^{2H-2} \ud s \ud t.
\end{split}
\end{equation*}
Therefore, the randomized periodogram is also a biased estimator of the quadratic variation $[X,X]$.
\end{itemize}

\begin{rmk} It would be interesting to know, whether the estimating based on `` discretized`` 
periodogram (or ''realized periodogram'') is less noisy than the realized quadratic variation estimator.
\end{rmk}
\begin{rmk}
 It is also interesting whether one can give an unbiased estimator of the quadratic variation of mixed Brownian-fractional Brownian motion.
\end{rmk}

\end{document}